\documentclass{article}
\usepackage{amsmath,amssymb}
 \newtheorem{thm}{Theorem}[section]
 \newtheorem{cor}[thm]{Corollary}

\begin{document}

 {\Large Supersolutions to degenerated logistic equation type}

\begin{center}
 
{Marcos Marv\'a.}

{Departamento de F\'isica y Matem\'aticas, Universidad de Alcal\'a,\\ 28871 Alcal\'a de Henares, Spain}

{marcos.marva@uah.es}
\end{center}


 


\begin{abstract}
In this work we provide a method for building up a strictly positive supersolution for the steady state of a degenerated logistic equation type, i.e., when the weight function vanishes on the boundary of the domain. This degenerated system is related in obtaining the so-called large solutions. Previously, this problem was handled as the limit case of non degenerated approaching problems.  Our method can be adapted straightforwardly to degenerated boundary value problems.

\end{abstract}


\section{Introduction}\label{intro}
In this work we show how to build up a positive strict supersolution to the boundary value problem
	\begin{equation}\label{pcm}
	\left\{\begin{array}{lll}
	-\Delta u=\lambda  m(x) u-a(x)f(x,u) & in  & \Omega,\\
	u=g(x) &  on & \partial \Omega,\\	
	\end{array}\right.
	\end{equation}
where $a\in C^{\mu}(\bar{\Omega}; [0,+\infty))$, $m\in C^{\mu}(\bar{\Omega}; \mathbb{R})$, $\mu\in(0,1)$, $a\equiv 0$ on $\partial \Omega$, $g\in C^{1+\mu}(\partial \Omega; [0,+\infty))$ is such that $g\geq 0$, $g\not\equiv 0$, 
$f\in C^{\mu,1+\mu}(\bar{\Omega}\times [0,+\infty); \mathbb{R})$, $f(x,0)=0$, $f(x,u)>0$, $\partial_{u}f(x,u)>0$ and $\lim_{k\nearrow+\infty}\frac{f(x,ku)}{k}=+\infty$ for each $u>0$ uniformly in $x\in\Omega$. The set $\Omega\subset\mathbb{R}^{N}$ is a bounded domain with  $\partial \Omega$ of class $C^{2}$.

Problem \eqref{pcm} (with $m\equiv 1$)  stands for the steady states of the generalized logistic growth law \cite{Mu93}. Typically, $u$ stands for the distribution in $\Omega$ of individuals of certain species, $\lambda$ stands for the net growth rate of $u$ modulated by $m(x)$, $a(x)$ simulates demography pressure and, 
along with $m(x)$, environmental heterogeneity. When $a\equiv 0$ somewhere, model \eqref{pcm} (with $m\equiv 1$)  becomes a Malthusian model in the region $\Omega_{0}:=\{x\in\Omega;\ a(x)=0\}$, known as {\it refuge}. We deal with the case $a\equiv 0$ on $\partial \Omega$, but it could be adapted to apply in the general case when $\Omega_{0}\neq\emptyset$ 
and $\Omega_{0}\subset\Omega$.\newline

In addition, solving problem \eqref{pcm} is an intermediate step to prove the existence of large solutions related to problem \eqref{pcm}, i.e., a function  $u\in C^{2,\nu}(\overline{\Omega})$ such that
	\begin{equation}\label{pcinfinito}
	\left\{\begin{array}{ll}
	-\Delta u=\lambda m(x) u-a(x) f(x,u)u &\,\, in\,\,\,\, \Omega , \\
	\displaystyle\lim_{x\rightarrow\partial \Omega}u(x)=\infty, & 	
	\end{array}
	\right.
	\end{equation}
see \cite{Lop00}. Large solutions come across from studies concerning combustion due to Keller \cite{Ke57} and Osserman \cite{Os57}. In fact, there are available uniqueness results for problem \eqref{pcinfinito} (
\cite{CR05}, \cite{DU04}, \cite{Lo04c} and references therein). As long as we know, problem \eqref{pcm} has been handled as the limit case of non degenerated approaching problems.
\newline

If we restrict ourselves in problem \eqref{pcm} to $a(x)\geq\gamma>0$ on $\overline{\Omega}$, it is easy to show that there exists $\epsilon>0$ small enough such that $\epsilon\varphi$ is a positive subsolution to problem \eqref{pcm}, being $\varphi$ the principal eigenfunction of $-\Delta$ in $\Omega$ under homogeneous Dirichlet boundary conditions. In this case a positive constant $K$ such that
	\begin{equation}\label{cotainf}K> \max \left\{ \underset{x\in\partial\Omega}{\max} \left\{g\right\},\ \dfrac{ \underset{x\in\overline{\Omega}}{\max}\left\{\lambda m(x)\right\}}{\underset{x\in\overline{\Omega}}{\min}\left\{a(x)\right\}}\right\}, 	
	\end{equation}
provides us with a supersolution to problem \eqref{pcm}. Now, enlarging $K$ if necessary,  we get an ordered pair $\left(\epsilon\varphi, K\right)$ of sub-supersolution. Thanks to a Theorem by Amann \cite{Am76}, there exists a solution to problem \eqref{pcm} between $\epsilon\varphi$ and $K$. Unfortunately, as soon as we let $a=0$ somewhere in $\bar{\Omega}$, we lose condition \eqref{cotainf}, so it is needed something  different from  a constant to get a supersolution to problem \eqref{pcm}. Different strategies have been used to avoid this problem: some authors approximate $\Omega$ by $\Omega_{n}:=\{x\in\Omega;\ \text{dist}(x,\partial \Omega)>1/n\}$ (for each $n$ $a(x)$ is uniformly bounded away from $0$ in $\Omega_{n}$, that is
 $a\geq \gamma>0$, see \cite{DH99} and \cite{Lop00}), and others approach $a(x)$ by means of $\frac{1}{n}+a(x)$, essentially in order to avoid condition $a\equiv 0$ on $\partial \Omega$ and to generate a sequence of solutions of approximate problems converging to a solution to problem \eqref{pcm}, see \cite{CR03a} and  \cite{GLS01}.

\section{Results}\label{results}
The result we present here allows us to obtain via sub and supersolution the existence of a solution to  \eqref{pcm}. Actually, the construction of the supersolution in $\Omega$ follows from an extension of the Faber \cite{Fab23} and  Krahn \cite{Kra25} inequality due to L\'opez-G\'omez \cite{Lop96} which provides us with a lower estimate for the main eigenvalue to $-\Delta$ on the domain $D$:
    \begin{equation}\label{EstimaAutovPPalJulian}\sigma_{1}[-\Delta,D]\geq
        \dfrac{\sigma_{1}[-\Delta,B_{1}(0)]\cdot
    |B_{1}(0)|^{2/N}}{|D|^{2/N}}, \end{equation}
 where $B_{1}(0):=\{x\in \mathbb{R}^{N};|x|\leq 1\}$, $\sigma_{1}[-\Delta,B_{1}(0)]$ is the principal eigenvalue of $-\Delta$ on $B_{1}(0)$ and $|D|$ stands for the Lebesgue's measure of $D$.
\begin{thm}
Consider the boundary value problem  \eqref{pcm}. Then, for each $\lambda\in\mathbb{R}$, there exists  an ordered pair $(\underline{u},\overline{u})$ consisting of a subsolution and a positive strict supersolution, where $\underline{u}(x)<\overline{u}(x)\ \forall x\in\Omega$ .
\end{thm}
\textbf{Proof}.-  Let us define
		$$\begin{array}{ccccc}
		  O_{\epsilon}:=\{x\in \mathbb{R}^{N};\ \text{dist}(x, \partial \Omega)< \epsilon\},
		  &  &  &  & \sigma_{\epsilon}:=\sigma[-\Delta,O_{\epsilon}],
		\end{array}
		$$
where, keeping in mind \eqref{EstimaAutovPPalJulian}, for a fixed $\lambda\in\mathbb{R}$ we can choose $\epsilon>0$ such that
    $$\underset{x\in\overline{\Omega}}{\max}\left\{\lambda m(x)\right\}<\sigma_{\epsilon}.$$
Let  $\varphi_{o}^{\epsilon}$ be the principal eigenfunction of $-\Delta$ under homogeneous boundary conditions on $O_{\epsilon}$.
 We define
    $$\Phi:= \left\{ \begin{array}{ll}
      \varphi_{o}^{\epsilon} & \text{in} \ \  \overline{\Omega}\cap O_{\epsilon/2},\\
      \Psi & \text{in} \ \Omega_{\epsilon/2}:=  \Omega\backslash
    (\overline{\Omega}\cap O_{\epsilon/2}),\\
    \end{array}\right.$$
where $\Psi\geq \tau>0$ is any regular enough function such that  $\Phi\in
C^{2}(\overline{\Omega})$.
We are looking for the existence of a constant $K>0$, large enough, such that function
    $$\overline{u}:=K\Phi $$
is a strict positive supersolution to problem \eqref{pcm}.\\

We estimate now $K$ on $\overline{\Omega}\cap O_{\epsilon/2}$. Keeping in mind that $ a(x)\geq 0$ and $K
\Phi=K\varphi_{o}^{\epsilon}$, it must be
$$\left\{
    \begin{array}{ll}
        K(-\Delta \varphi_{o}^{\epsilon})\geq \lambda m K\varphi_{o}^{\epsilon}
             -a f(x,K\varphi_{o}^{\epsilon})&  \text{in}\ \ \overline{\Omega}\cap O_{\epsilon/2},\\
        K\varphi_{o}^{\epsilon}\geq g>0 & \text{on} \ \ \partial  \overline{\Omega}\cap O_{\epsilon/2}=\partial   \Omega.\\
        \end{array}
    \right.$$
In the interior of the domain it should happen
    $$ K(-\Delta \varphi_{o}^{\epsilon})\geq \lambda m K\varphi_{o}^{\epsilon}
          -a f(x,K\varphi_{o}^{\epsilon}),$$
which is equivalent to
        $$ K\sigma_{\epsilon} \varphi_{o}^{\epsilon}\geq \lambda m K\varphi_{o}^{\epsilon}
             -a f(x,K\varphi_{o}^{\epsilon}) $$
by the definition of $\varphi_{o}^{\epsilon}$. Rearranging terms we get
    \begin{equation}\label{K1} af(x,K\varphi_{o}^{\epsilon})\geq K\varphi_{o}^{\epsilon}\left(\lambda m- \sigma_{\epsilon}\right),\end{equation}
 where the left hand side is non negative and the right one is strictly negative because of the selection of $\epsilon$. Therefore, \eqref{K1} is a proper inequality $\forall K>0$. 
By construction, on the boundary we have
$\left.\varphi_{o}^{\epsilon}\right|_{\partial \Omega},g>0$ and therefore, whenever
         \begin{equation}\label{K2}K>\max_{x\in\partial \Omega}\left\{\dfrac{g(x)}{\varphi_{o}^{\epsilon}(x)}\right\},\end{equation}
the boundary condition is satisfied (with strict inequality).\\

In the domain $\Omega_{\epsilon/2}$ it should be verified
	\begin{equation}\label{supersolint}K(-\Delta \Psi)\geq \lambda m K\Psi-a f(x,K\Psi).\end{equation}
rearranging terms, expression \eqref{supersolint} is equivalent to
	\begin{equation}\label{K3} \dfrac{f(x,K\Psi)}{K}\geq \underset{x\in\overline{\Omega}_{\epsilon/2}}{\max}\left\{\dfrac{\Delta\Psi+\lambda m\Psi}{a}\right\}.\end{equation}
As $a(x)>0\ \ \forall x\in\Omega_{\epsilon/2}$, function $\frac{\Delta\Psi+\lambda m\Psi}{a}$ reaches its maximum (provided a regular $\Psi$). By hypothesis on $f$, there exists $K_{0}$ such that $\forall K>K_{0}$ condition \eqref{K3} holds. Thus, there exists $K>0$ such that function $K\Phi$ is a strict positive supersolution to \eqref{pcm}. Function $\underline{u}:=0$ is a subsolution to problem \eqref{pcm}, and $\underline{u}(x)<\overline{u}(x)\  \forall x\in\bar{\Omega}$. $\blacksquare$

\begin{cor} For each $\lambda\in\mathbb{R}$, problem \eqref{pcm} has an unique positive solution.
\end{cor}
\textbf{Proof}.- Since we have an ordered pair formed by a subsolution and a supersolution, a Theorem due to Amann \cite{Am76} guaranties the existence of a  solution $u$ to problem \eqref{pcm} such that
    $0\leq u\leq K\Psi.$
Uniqueness follows, for instance, from \cite{Lop00}.   $\blacksquare$

\textbf{Acknowledgement}.- The author thanks E. G\'omez-Hidalgo  for her fruitful comments and support. The author is  partially supported by Ministerio de Ciencia e Innovaci\'on (Spain), 
projects MTM2011-24321 and MTM2011-25238.

\end{document}